\documentclass[11pt]{article}
\usepackage{amsmath}
\usepackage{amssymb}
\usepackage{amscd}
\newtheorem{theorem}{Theorem}[section]
\newtheorem{proposition}[theorem]{Proposition}
\newtheorem{corollary}[theorem]{Corollary}
\newtheorem{lemma}[theorem]{Lemma}

\newtheorem{remark}[theorem]{Remark}
\newtheorem{remarks}[theorem]{Remarks}
\newenvironment{proof}{\begin{trivlist} \item[]{\bf Proof.}}
{\par\hfill $\square$\end{trivlist}}
\newcommand{\DSH}{{\rm DSH}}

\newcommand{\Ball}{{\rm B}}

\newcommand{\C}{\mathbb{C}}

\newcommand{\R}{\mathbb{R}}

\newcommand{\h}{{\rm h}}
\newcommand{\M}{{\cal M}}
\renewcommand{\L}{{\rm L}}

\newcommand{\Cr}{{\rm C}}

\newcommand{\dc}{{\rm d^c}}
\newcommand{\ddc}{{\rm dd^c}}
\renewcommand{\d}{{\rm d}}

\newcommand{\see}{{see }}
\newcommand{\ie}{{i.e. }}
\newcommand{\lov}{{\rm lov}}
\newcommand{\vol}{{\rm vol}}

\newcommand{\Lone}{{{\rm L}^1}}
\newcommand{\Ltwo}{{{\rm L}^2}}
\newcommand{\Linfty}{{{\rm L}^\infty}}

\newcommand{\dist}{{\rm dist}}

\newcommand{\supp}{{\rm supp}}

\title{Regularization of currents and entropy}
\author{Tien-Cuong Dinh and Nessim Sibony}
\date{}
\begin{document}
\maketitle
\small
\noindent
{\bf Abstract.} Let $T$ be a positive closed 
$(p,p)$-current on a compact K\"ahler manifold $X$. Then,
there exist smooth positive closed 
$(p,p)$-forms $T_n^+$ and $T_n^-$ such that $T^+_n-T^-_n\rightarrow T$
weakly. Moreover, $\|T^\pm_n\|\leq c_X\|T\|$ where $c_X>0$ is a constant 
independent of $T$.  We also extend this result to positive pluriharmonic currents.
Then we study the wedge product of positive closed $(1,1)$-currents having
continuous potential with positive pluriharmonic currents.
As an application, we 
give an estimate for the topological entropy of meromorphic maps 
on compact K\"ahler manifolds. 
\\
\ 
\\
{\bf R\'esum\'e.} Soit $T$ un $(p,p)$-courant positif ferm\'e 
sur une vari\'et\'e k\"ahl\'erienne compacte $X$. Alors, il existe 
des $(p,p)$-formes
lisses, positives ferm\'ees $T_n^+$ et $T_n^-$
telles que $T_n^+-T_n^-\rightarrow T$ faiblement. De plus, on a 
$\|T^\pm_n\|\leq c_X\|T\|$ o\`u $c_X>0$ est une constante ind\'ependante de $T$.
Nous montrons aussi ce r\'esultat pour les courants positifs pluriharmoniques.
Nous \'etudions \'egalement le produit ext\'erieur de $(1,1)$-courants 
positifs ferm\'es \`a potentiel continu avec des courants pluriharmoniques positifs.
Comme application, nous donnons une estimation de l'entropie topologique 
des applications m\'eromorphes d'une vari\'et\'e k\"ahl\'erienne compacte.

\normalsize
\section{Introduction}
Let $(X,\omega)$ be a compact K\"ahler manifold of dimension $k$.
Demailly \cite{Demailly2} has shown that for a positive closed
$(1,1)$-current $T$ on $X$, there exist smooth positive closed
$(1,1)$-forms $T_n^+$ which converge weakly (\ie in the sense
of currents) to $T+c\omega$
where $c>0$ is a constant.
Moreover, there is a constant $c_X>0$, independent of $T$, such
that $\|T^+_n\|$ and $c$ are bounded by $c_X\|T\|$. 
We refer to Demailly's papers \cite{Demailly1, Demailly2} for the basics 
on currents on complex manifolds. Recall that the mass
of a positive $(p,p)$-current $S$ is defined 
by $\|S\|:=\int_X S\wedge \omega^{k-p}$. 
Our main result is the following theorem where 
the positivity can be understood in the
weak or strong sense.
\begin{theorem} Let $(X,\omega)$ be a compact K\"ahler manifold
of dimension $k$. Then, for every
positive closed $(p,p)$-current $T$ on $X$, there exist smooth closed
$(p,p)$-forms $T^+_n$ and $T^-_n$
such that $T^+_n-T^-_n$ converge weakly
to the current $T$. Moreover, $\|T^\pm_n\|\leq
c_X\|T\|$ where $c_X>0$ is a constant independent of $T$. 
\end{theorem} 
We deduce from this theorem the following corollary which is proved 
in \cite{DinhSibony1} 
for projective manifolds.
\begin{corollary} Let $(X,\omega)$ be a compact K\"ahler manifold
of dimension $k$. Then, for every
positive closed $(p,p)$-current $T$ on $X$, there exist smooth closed
$(p,p)$-forms $T^+_n$
which converge weakly to a current $T'$
with $T'\geq T$. Moreover, $\|T^+_n\|\leq
c_X\|T\|$ and $\|T'\|\leq c_X\|T\|$
where $c_X>0$ is a constant independent of $T$.
\end{corollary}
\par
Let $(X',\omega')$ be  another compact K\"ahler manifold of
dimension $k'\geq k$ and let  $\Pi:X'\longrightarrow X$ be a surjective
holomorphic map. We want to define the
pull-back of the current $T$ by the map $\Pi$. 
When $\Pi$ is a finite map, this problem is studied in \cite{Meo, DinhSibony3}.
In general, the map $\Pi$ is a submersion only in the complement of an analytic subset
$\Cr$ of $X'$. Let $\pi$ denote the restriction of $\Pi$ to $X'\setminus
\Cr$. 
Then, $\pi^*(T)$ is well defined and 
is a positive closed $(p,p)$-current on $X'\setminus \Cr$. Let $(T_n^+)$ 
and $c_X$ be
as in Corollary 1.2. Define $S_n:=\Pi^*(T_n^+)$. The $(p,p)$-forms 
$S_n$ are smooth and positive 
on $X'$. Their classes in $H^{p,p}(X',\C)$ are bounded since 
$(\|T^+_n\|)$ is bounded. It follows
that $(\|S_n\|)$ is bounded.
Taking a subsequence, we can assume that $S_n$ 
converge to a current $S$. We also have $S\geq\pi^*(T)$ on $X'\setminus \Cr$.
In particular, $\pi^*(T)$ has finite mass.
Following Skoda \cite{Skoda}, 
the trivial extension $\widetilde{\pi^*(T)}$ of $\pi^*(T)$ on $X'$ is a positive
closed current. So, we have the following corollary.
\begin{corollary} Let $X$, $X'$, $\Pi$, $\pi$ and $T$ be as above. Then, 
the positive 
current $\widetilde{\pi^*(T)}$ is well defined and closed. Moreover, there
exists a constant $c_\Pi>0$ independent of $T$ such that 
$\|\widetilde{\pi^*(T)}\|\leq c_\Pi\|T\|$. The map $T\mapsto \widetilde{\pi^*(T)}$
is l.s.c. in the sense that if $T_n\rightarrow T$, then any cluster point $\tau$ of
$\big(\widetilde{\pi^*(T_n)}\big)$ satisfies $\tau\geq\widetilde{\pi^*(T)}$. 
\end{corollary}
In \cite{Meo}, M\'eo gave an example which shows that, 
in general, when $X$ and $X'$ are
not compact, the current $\pi^*(T)$ on $X\setminus\Cr$ is not always 
of bounded mass near $\Cr$.
\par
Consider a dominating meromorphic self-map $f:X\longrightarrow X$ of $X$. Define 
$f^n:=f\circ\cdots\circ f$ ($n$ times) the $n$-th iterate of $f$. 
We refer to the survey \cite{Sibony} for the theory of iteration 
of meromorphic maps.
Let $I_n$ be the 
indeterminacy set of $f^n$. Then $I_n$ is an analytic subset of 
codimension $\geq 2$ of 
$X$. 
Denote by $\Omega_f$ the set of points $x\in X\setminus I_1$ such that
$f^n(x)\not\in I_1$ for every $n\geq 1$. A subset $F\subset \Omega_f$ is called
{\it $(n,\epsilon)$-separeted}, $\epsilon>0$, if
$$\max_{0\leq i\leq n-1}\dist(f^i(x),f^i(y))\geq\epsilon \mbox{ for }
x,y\in F \mbox{ distinct.}$$
The {\it topological entropy} $\h(f)$ (\see \cite{Bowen}) is defined by
$$\h(f):=\sup_{\epsilon>0}\left(\limsup_{n\rightarrow\infty} \frac{1}{n} 
\log\max \big\{\#F,\ F \ (n,\epsilon)\mbox{-separated}\big\} \right).$$
Let $\Gamma_n$ be the closure in $X^n$ of the set of points     
$$(x,f(x),\ldots, f^{n-1}(x)),\ \ x\in\Omega_f.$$
This is an analytic subset of dimension $k$ of $X^n$. Let $\Pi_i$ be the 
canonical projections of $X^n$ on its factors. We consider
on $X^n$ the K\"ahler metric $\omega_n:=\sum\Pi_i^*(\omega)$.
Define following Gromov \cite{Gromov2},
\begin{eqnarray}
\lov(f) & := & \limsup_{n\rightarrow\infty}\frac{1}{n}\log(\vol(\Gamma_n))
=\limsup_{n\rightarrow\infty}\frac{1}{n}
\log\left(\int_{\Gamma_n} \omega_n^k\right).
\end{eqnarray}
Define also {\it the dynamical degree of order $p$} of $f$ by 
\begin{eqnarray}
d_p & := & \limsup_{n\rightarrow\infty}\left(\int_{X\setminus I_n} 
f^{n*}(\omega^p)
\wedge \omega^{k-p}\right)^{1/n}.
\end{eqnarray}
Using an inequality of Lelong \cite{Lelong}, Gromov \cite{Gromov2} proved that 
$\h(f)\leq \lov(f)$. Following Gromov and Yomdin \cite{Yomdin,Gromov1, Gromov2}, 
we have
$$\h(f)=\lov(f)=\max_{1\leq p\leq k} \log d_p$$ 
when  $f$ is a holomorphic map.
Using Corollary 1.2, we prove, in the same way as in \cite{DinhSibony1}, that 
the sequences in (1)(2) are convergent and that the dynamical degrees 
$d_p$ are bimeromorphic
invariants of $f$. More precisely, if $\Pi:X'\longrightarrow X$ is a 
bimeromorphic map between compact K\"ahler manifolds, the dynamical degrees of 
$\Pi^{-1}\circ f\circ\Pi$ are equal to $d_p$.
Using Corollary 1.2, we also get the
following result.
\begin{theorem} Let $f$ be a dominating meromorphic self-map on a compact K\"ahler manifold
$X$ of dimension $k$. Let $d_p$ be the dynamical degrees of $f$. Then
$$\h(f)\leq\lov(f)=\max_{1\leq p\leq k} \log d_p.$$
\end{theorem}
This theorem gives a partial answer to a conjecture of Friedland \cite{Friedland}
which says that $\h(f) =\max_{1\leq p\leq k} \log d_p$.
Theorem 1.4 is already proved in \cite{DinhSibony1} for rational maps 
on projective
manifolds.
Corollary 1.2 permits to extend the proof
to the case of compact K\"ahler manifolds. One can also extend
some results on 
meromorphic correspondences or transformations, which are proved 
in the projective case in \cite{DinhSibony2}
(\see also \cite{Dinh}). 

In the last two sections, we extend Theorem 1.1 to positive 
pluriharmonic currents
and currents of class DSH. 
We also study the intersection of such  currents with 
positive closed $(1,1)$-currents.

We thank the referee for his constructive observations that helped to improve
the exposition.
\section{A classical lemma}
We will give here a classical lemma that we use in Section 3. 
Let $\Ball$ denote the unit ball in $\R^m$. Let 
$K(x,y)$ be a function with compact support in $\Ball\times\Ball$,
smooth in $\Ball\times\Ball\setminus\Delta$ where $\Delta$ 
is the diagonal 
of $\Ball\times\Ball$. Assume that, for every $(x,y)$
\begin{eqnarray}
|K(x,y)| & \leq &  A|x-y|^{2-m}
\end{eqnarray} 
where $A>0$ is a constant and $x=(x_1,\ldots,x_m)$ are coordinates of $\R^m$. 
Observe that for every $y$
\begin{eqnarray}
\|K(.,y)\|_{\L^{1+\delta}} & \leq &  A'
\end{eqnarray} 
for some $\delta>0$ and $A'>0$.
Assume also that for every $x$, $y$ 
\begin{eqnarray}
|\nabla K(x,y)| & \leq &  A|x-y|^{1-m}.
\end{eqnarray}

In this section, we identify $\nu$, a current of degree $0$ and of order $0$,
with the current of degree $m$,  $\nu\d y_1\wedge\ldots\wedge \d y_m$.
Let $\M$ denote the set of Radon measures on $\R^m$.
We define a linear operator $P$ on $\M$ 
by: 
$$P\mu(x):=\int_{y\in\R^m} K(x,y)\d\mu(y).$$
Observe that the function 
$P\mu$ has support in
$\Ball$. We have the following lemma.
\begin{lemma} The operator $P$ maps continuously 
$\M$ into ${\rm L}^{1+\delta}$. It also maps continuously 
$\L^p$ into $\L^{q}$,  
$\Linfty$ into ${\cal C}^0$ and ${\cal C}^0$ 
into ${\cal C}^1$,
where $q=\infty$ if $p^{-1}+(1+\delta)^{-1}\leq 1$ and
$p^{-1}+(1+\delta)^{-1}= 1 + q^{-1}$ otherwise. 
\end{lemma}
All the assertions are easy to deduce from (3)(4)(5) and the 
H\"older inequality.
\section{Proof of Theorem 1.1}
Let $\Delta$ denote the diagonal of $X\times X$. We first give a {\it weak 
regularization} of the current of integration $[\Delta]$.
Let $\widetilde{X\times X}$ denote the blow-up of $X\times X$ along $\Delta$.
Following Blanchard \cite{Blanchard}, $\widetilde{X\times X}$ is a K\"ahler 
manifold. 
Let $\pi: \widetilde{X\times X}\longrightarrow X\times X$ be the canonical 
projection and $\widetilde{\Delta}:=\pi^{-1}(\Delta)$.
Then $\widetilde{\Delta}$ is a smooth 
hypersurface in $\widetilde{X\times X}$. 
If $\gamma$ is a closed strictly positive $(k-1,k-1)$-form 
on $\widetilde{X\times X}$, then $\pi_*(\gamma\wedge [\widetilde{\Delta}])$ 
is a non-zero positive closed $(k,k)$-current on $X\times X$ supported on $\Delta$.
So, it is a multiple of $[\Delta]$. We choose $\gamma$ so that 
$\pi_*(\gamma\wedge [\widetilde{\Delta}])=[\Delta]$. We will use the following regularization of $[\widetilde
\Delta]$.    
\par
Since $[\widetilde\Delta]$ is a positive closed $(1,1)$-current, there exist
a quasi-p.s.h. function $\varphi$ and a smooth closed $(1,1)$-form $\Theta'$ such that
$\ddc\varphi=[\widetilde\Delta]-\Theta'$. Recall that 
$\dc:=\frac{i}{2\pi}(\overline\partial-\partial)$.
Demailly's regularization theorem \cite{Demailly2} 
implies the existence of
smooth functions $\varphi_n$ and of a smooth positive closed 
$(1,1)$-form $\Theta$ 
on $\widetilde{X\times X}$ such that
\begin{enumerate}
\item[$\bullet$] $\ddc\varphi_n\geq -\Theta$;
\item[$\bullet$] $\varphi_n$ decrease to $\varphi$.
\end{enumerate}

In this case, independently of Demailly's theorem, we can construct the functions $\varphi_n$ as follows.
Observe that $\varphi$ is smooth out of $\widetilde\Delta$ and 
$\varphi^{-1}(-\infty)=\widetilde\Delta$.
Let $\chi:\R\cup\{-\infty\}\rightarrow \R$
be a smooth increasing convex function such that $\chi(x)=0$ on $[-\infty,-1]$,
$\chi(x)=x$ on $[1,+\infty[$ and $0\leq\chi'\leq 1$. Define $\chi_n(x):=\chi(x+n)-n$
and $\varphi_n:=\chi_n\circ\varphi$. The functions $\varphi_n$ 
are smooth decreasing to $\varphi$ and we have
\begin{eqnarray}
\ddc\varphi_n & = & (\chi''_n\circ\varphi)\d\varphi\wedge \dc\varphi +
(\chi_n'\circ\varphi) \ddc\varphi\nonumber \\
& \geq &  (\chi_n'\circ\varphi) \ddc\varphi = -(\chi_n'\circ \varphi)\Theta'\geq -\Theta
\end{eqnarray}
where we choose the smooth positive closed form $\Theta$ big enough so that $\Theta-\Theta'$ is positive.

Define $\Theta^+_n:=\ddc\varphi_n+\Theta$ and $\Theta_n^-:=\Theta-\Theta'$ then
$\Theta^+_n-\Theta^-_n\rightarrow [\widetilde{\Delta}]$. 
We have
$\|\Theta^\pm_n\|\leq c_0$ where $c_0>0$ is a constant.
The forms $\Theta_n^\pm$ are smooth.
Define 
$$\widetilde{K}^\pm_n:=\gamma\wedge \Theta^\pm_n
\ \mbox{ and }\ K^\pm_n:=\pi_*(\widetilde{K}^\pm_n).$$
The $(k,k)$-forms $K_n^\pm$ are positive closed with coefficients in $\Lone$ and 
smooth out of $\Delta$.
We also have $K^+_n-K^-_n\rightarrow [\Delta]$ weakly and 
$\|K^\pm_n\|\leq c_1$, $c_1>0$. This is what we call
a {\it weak regularization} of $[\Delta]$. 
We will use $K^\pm_n$ to regularize the current $T$. 
The following lemma shows that the coefficients of $K_n^\pm$ satisfy inequalities of type (3) and (5) for $m=2k$.
Then, the singularities of $K_n^\pm$ are the same  
than for the Bochner-Martinelli kernel.

\begin{lemma} Let $(x,y)=(x_1,\ldots,x_k,y_1,\ldots ,y_k)$, $|x_i|<3$, $|y_i|<3$, 
be local holomorphic coordinates of a chart of $X\times X$ such that
$\Delta=(y=0)$ in that chart. Let $H_n^\pm$ be a coefficient of $K_n^\pm$ in these
coordinates. Then, there exists a constant  $A_n>0$, depending on $n$, such that 
$$|H_n^\pm(x,y)|\leq A_n|y|^{2-2k} \mbox{\ \ and \ } |\nabla H_n^\pm|\leq A_n|y|^{1-2k}$$ 
for $|x_i|\leq 1$, $|y_i|\leq 1$ and $y\not=0$. 
\end{lemma}
\begin{proof}
By symmetry, it is sufficient to consider $(x,y)$ in the open sector $S$ defined by the inequalities 
$|x_i|<3$, $|y_i|<3$, $|y_i|< 3|y_1|$ and prove the estimates in the sector $S'$ defined by
$|x_i|<2$, $|y_i|<2$ and $|y_i|< 2|y_1|$ (we can assume that $y_1$ is the largest coordinate of the point $y\not=0$).
Let $\widetilde S$ and $\widetilde S'$ be the interiors of $\pi^{-1}(\overline S)$ 
and of $\pi^{-1}(\overline S')$ respectively.
We consider the coordinate system $(x,Y)$ of $\widetilde S$ with $|x_i|<3$, $Y_1=y_1$ and $Y_i=y_i/y_1$,
$|y_1|<3$, $|y_i|<3|y_1|$ for $i=2,\ldots,k$. We have $\pi(x,Y)=(x,y)$ for $(x,y)\in S$.
The equation of $\widetilde\Delta$ in $\widetilde S$ is $Y_1=0$.

Since $\widetilde K_n^\pm$ are smooth on $\widetilde S$, they are finite sums of forms of type 
$$\Phi(x,Y)=L(x,Y)\d x_I \wedge \d \overline x_{I'} \wedge \d Y_J \wedge \d \overline Y_{J'}$$
where $L$ is a smooth function, $I$, $I'$, $J$, $J'$ are subsequences of 
$\{1,\ldots,k\}$ and $\d x_I=\d x_{i_1}\wedge \ldots \wedge \d x_{i_m}$
if $I=\{i_1,\ldots, i_m\}$. Hence, in $S$ the forms $K_n^\pm$ are finite sums of forms of 
type $\pi_*(\Phi)$. 

Observe that $\pi_*(\Phi)$ is obtained from $\Phi(x,Y)$ 
replacing $Y_1$ by $y_1$ and $Y_i$ by $y_i/y_1$.
There are here at most $2k-2$ factors of the form $\d(y_i/y_1)=\d y_i/y_1-y_i\d y_1/y_1^2$ or their conjugate. 
Hence, the coefficients of $\pi_*(\Phi)$ on 
$S$ are finite sums of  
$$L(x,y_1,y_2/y_1,\ldots, y_k/y_1)P(y)y_1^{-m}\overline y_1^{-n}$$
where $P$ is a homogeneous polynomial such that $\deg(P)+2k-2\geq m+n$. 
Since $\widetilde S'\Subset \widetilde S$, 
$L$ is bounded on $\widetilde S'$ and $L(x,y_1,y_2/y_1,\ldots, y_k/y_1)$ is 
bounded on $S'$. The first estimate of the lemma follows.

For the second estimate, it is sufficient to observe that the coefficients in the gradient of
$$L(x,y_1,y_2/y_1,\ldots, y_k/y_1)P(y)y_1^{-m}\overline y_1^{-n}$$
are combinations of functions of the same type with homogeneous polynomials $P$ such that $\deg(P)+2k-1\geq m+n$.
\end{proof}

Define
\begin{eqnarray}
T^\pm_n(x) & := & \int_{y\in X} K^\pm_n(x,y)\wedge T(y).
\end{eqnarray}
Let $\pi_i$ denote the canonical projections of $X\times X$ on its factors.
We have 
\begin{eqnarray}
T^\pm_n & := & (\pi_1)_* \big(K^\pm_n\wedge \pi_2^*(T)\big).
\end{eqnarray}
Observe that $\pi_2^*(T)$ is well defined since $\pi_2$ is a submersion. The
currents
$K^\pm_n\wedge \pi_2^*(T)$ are positive closed and well defined on 
$X\times X\setminus\Delta$. They are of finite mass since, 
for each $n$, $\|K_n^+(.,y)\|_\Lone$ is uniformly bounded. 
A priori, the mass depends on $n$.
By Skoda's extension theorem 
\cite{Skoda}, their trivial extensions are positive and closed. It follows that 
$T^\pm_n$ are well defined and are positive closed currents on $X$.
The use of Skoda theorem can be replaced by an argument similar to the one 
in the proof of the following 
lemma. 
\begin{lemma} The currents $T_n^+-T_n^-$ converge weakly to $T$ 
when $n\rightarrow\infty$. Moreover, $\|T_n^\pm\|\leq c\|T\|$ where $c>0$ is a constant 
independent of $n$ and $T$.
\end{lemma}
\begin{proof} Define $\Pi:=\pi_2\circ\pi$. Observe that $\Pi$ is a submersion
from $\widetilde{X\times X}$ onto $X$ and 
$\Pi_{|\widetilde\Delta}$ is a submersion from $\widetilde\Delta$
onto $X$.   
Indeed, consider charts $U\Subset V'\subset X$ that we identify with open sets in $\C^k$.
Assume that $U$ is small enough and $0\in U$.
We can, using the change of coordinates $(z,w)\mapsto (z-w,w)$ on $V'\times U$, 
reduce to the
product situation $V\times U$, $U\Subset V\subset \C^k$
where $\Delta$ is identified to $\{0\}\times U$.
The blow-up along $\{0\}\times U$ is still a product. So $\Pi^*$ of a current is just
integration on fibers. We can use this local model for the assertions below.

The potential of $\widetilde{\Delta}$ is integrable with respect to $\Pi^*(T)$ since its singularity
is like $\log\dist(z,\widetilde{\Delta})$ and this function has bounded integral
on fibers of $\Pi$. In particular, $[\widetilde{\Delta}]\wedge\Pi^*(T)$ is well
defined and is equal to $(\Pi_{|\widetilde\Delta})^*(T)$, and 
$[\widetilde{\Delta}]$ has no mass for $\Pi^*(T)$ nor for $\widetilde{K}_n^\pm
\wedge\Pi^*(T)$. We then have
\begin{eqnarray}
K^\pm_n\wedge \pi_2^*(T)=\pi_*(\widetilde{K}^\pm_n\wedge \Pi^*(T))
\end{eqnarray}
since the formula is valid out of $\Delta$.
The potentials of $\widetilde{K}_n^+$ are decreasing and 
the currents $\widetilde{K}_n^-$
are independent of $n$,
hence
\begin{eqnarray}
\widetilde{K}^+_n\wedge \Pi^*(T) - \widetilde{K}^-_n\wedge \Pi^*(T)
\rightarrow \gamma\wedge [\widetilde\Delta]\wedge \Pi^*(T)=
\gamma\wedge (\Pi_{|\widetilde\Delta})^*(T).
\end{eqnarray}
Since $\pi_{|\widetilde{\Delta}}$ is a submersion onto $\Delta$,  
we have $(\Pi_{|\widetilde\Delta})^*(T)=(\pi_{|\widetilde\Delta})^*(\pi_{2|\widetilde\Delta})^*(T)$. Hence
$$\pi_*\big(\gamma\wedge(\Pi_{|\widetilde\Delta})^*(T)\big)=(\pi_{2|\Delta})^*(T).$$
This and (9) (10) imply that 
$$K^+_n\wedge \pi_2^*(T) - K^-_n\wedge \pi_2^*(T)
\rightarrow (\pi_{2|\Delta})^*(T).$$
Taking the direct image under $\pi_1$ gives $T_n^+-T_n^-\rightarrow T$.

Since $\Pi$ is a submersion, $\|\Pi^*(T)\|\leq c_2\|T\|$ where $c_2>0$ is 
independent of $T$. Observe that since $\widetilde K_n^\pm$ are smooth we can 
compute $\|\widetilde K_n^\pm\wedge \Pi^*(T)\|$ cohomologically. The cohomological
classes of $\widetilde K_n^\pm$ are bounded, hence there exists a constant $c_3>0$ 
such that $\|\widetilde K_n^\pm\wedge \Pi^*(T)\|\leq c_3\|T\|$. 
It follows that
$$\|T^\pm_n\|=\|(\pi_1)_*\pi_*(\widetilde K_n^\pm\wedge \Pi^*(T))\|\leq c\|T\|$$
where $c>0$ is independent of $n$ and $T$. 
\end{proof}

The proof of Theorem 1.1 is completed 
by the following three steps.
\begin{trivlist}
\item[]{\bf Step 1.} We show first that we can choose in 
Theorem 1.1 forms $T^\pm_n$ with $\Lone$ coefficients. 
Define $T^\pm_n$ as in (7)(8). 
We use partitions of unity of $X$ and of $X\times X$ in order to reduce 
the problem
to the case of $\R^m$. 
Following Lemmas 2.1 and 3.1, the forms $T^\pm_n$ have $\Lone$ coefficients. 
Lemma 3.2 implies that $T^+_n-T^-_n\rightarrow T$ and
$\|T^\pm_n\|\leq c\|T\|$.
Of course, in general, $T_n^+-T_n^-$ do not converge in $\Lone$ since 
the constants $A_n$ in Lemma 3.1 depend on $n$.

\item[] {\bf Step 2.} We can now assume that $T$ is a form with $\Lone$ 
coefficients. Define $T^\pm_n$ as in (7)(8). Lemmas 2.1 and 3.1 imply that 
$T^\pm_n$ are forms with coefficients in $\L^{1+\delta}$.
We also have $T^+_n-T^-_n\rightarrow T$ and $\|T_n^\pm\|\leq c\|T\|$.
Hence, we can assume that $T$ is a form with $\L^{1+\delta}$ coefficients.
We repeat this process $N$ times with $N\geq \delta^{-1}$. Lemmas 2.1, 3.1 and 3.2 allow
to reduce the problem to the case where $T$ is a form with $\Linfty$ coefficients.
If we repeat this process two more times, we can assume that $T$ is a ${\cal C}^1$
form.
\item[] {\bf Step 3.} Now assume that $T$ is of class ${\cal C}^1$.
We can also assume that $T$ is strictly positive. 
Let $\Omega$ be a smooth real closed $(p,p)$-form cohomologous to $T$. 
Using standard Hodge theory \cite{Demailly1}, there is a
real $(p-1,p-1)$-form $u$ of class ${\cal C}^2$ such that $T=\Omega+\ddc u$. 
Let $(u_n)$ be a sequence of real smooth $(p-1,p-1)$-forms such that 
$u_n\rightarrow u$ in ${\cal C}^2$ norm. The current $T_n:=\Omega +\ddc u_n$ 
converges to $T$ in ${\cal C}^0$ norm. Moreover, 
$T_n$ is positive for $n$ big enough since $T$ is strictly  
positive. This completes the proof of Theorem 1.1.
\end{trivlist}
\hfill $\square$
\section{Pluriharmonic currents}

In this section, we extend Theorem 1.1 to positive pluriharmonic currents, \ie 
positive $\ddc$-closed currents. We have 
the following result which is new even for bidegree $(1,1)$ currents.

\begin{theorem} Let $T$ be a positive $\ddc$-closed $(p,p)$-current 
on a compact K\"ahler 
manifold $(X,\omega)$. Then there exist smooth positive $\ddc$-closed forms 
$T_n^\pm$ such that $T_n^+-T_n^-\rightarrow T$. Moreover, $\|T_n^\pm\|\leq c_X\|T\|$
where $c_X>0$ is a constant independent of $T$.  
\end{theorem}
We deduce from this theorem the following corollary.

\begin{corollary} Let $X$, $X'$, $\Pi$, $\pi$ and $\Cr$ be as in Corollary 1.3. 
If $T$ is as in Theorem 4.1, then the positive $\ddc$-closed current 
$\widetilde{\pi^*(T)}$ is well defined. Moreover the operator $T\mapsto 
\widetilde{\pi^*(T)}$ is l.s.c. and $\|\widetilde{\pi^*(T)}\|\leq c_\Pi\|T\|$
where $c_\Pi>0$ is a constant independent of $T$.
\end{corollary}
To prove the corollary, 
observe that by Theorem 4.1, the positive pluriharmonic current 
$\pi^*(T)$, which is well defined on $X'\setminus \Cr$,
has finite mass. Following  Alessandrini-Bassanelli \cite{AlessandriniBassanelli}, 
$\widetilde{\pi^*(T)}$ satisfies 
$\ddc \widetilde{\pi^*(T)}\leq 0$. Then, Stokes Theorem implies that 
$\ddc \widetilde{\pi^*(T)}=0$. 
\par
\
\\
{\it Proof of Theorem 4.1.} We use the same idea as in Section 3. 
Clearly $T_n^\pm$ given by (7)(8)(9) are pluriharmonic positive currents.
We only need to check 
that $T_n^+-T_n^-\rightarrow T$. The rest of proof is the same as in Theorem 1.1.

Let $\varphi$ and $\varphi_n$ be  q.p.s.h. functions as in Section 3.
We want to prove the analog of (10):
\begin{eqnarray}
(\ddc\varphi_n+\Theta')\wedge\Pi^*(T)\rightarrow (\Pi_{|\widetilde\Delta})^*(T)
\end{eqnarray}
The problem is local. Define $S:=(\Pi_{|\widetilde\Delta})^*(T)$. We choose as in Lemmas 3.1 and 3.2 local 
holomorphic coordinates $(x_1,\ldots,x_{2k})$ of an open set 
$W\subset\widetilde{X\times X}$, $|x_i|< 1$, 
so that in $W$
\begin{enumerate}
\item[$\bullet$] $\widetilde\Delta=\{x_{2k}=0\}$; hence $\psi:=\varphi-\log|x_{2k}|$
is smooth and $\ddc\psi=-\Theta'$;
\item[$\bullet$] $\Pi(x_1,\ldots,x_{2k})=(x_1,\ldots,x_k)$.
\end{enumerate}
Define $\tau(x_1,\ldots,x_{2k}):=(x_1,\ldots,x_{2k-1})$. 
Since $\Pi=\Pi_{|\widetilde\Delta}\circ\tau$, we have $\Pi^*(T)=\tau^*(S)$ in $W$.

Observe that $(\ddc \varphi_n+\Theta')\wedge \tau^*(S)$ is supported in $(\varphi<-n+2)$ and, by (6),
$(\ddc\varphi_n+\Theta')\wedge \tau^*(S)\geq 
(1-\chi_n'\circ\varphi)\Theta'\wedge \tau^*(S)$. The definition of $\chi_n$
implies that the measures $(1-\chi_n'\circ\varphi)\Theta'\wedge \tau^*(S)$ tend to 0. Hence,
every limit value of $(\ddc\varphi_n+\Theta')\wedge\tau^*(S)$ is 
a positive $\ddc$-closed current supported in $\widetilde\Delta$. 
Following Bassanelli
\cite{Bassanelli}, it is a current on $\widetilde\Delta$ 
(this is true for every positive current $T$ supported in $\widetilde\Delta$
such that $\ddc T$ is of order 0).
Hence, in order 
to prove (11) we only have to check that
$$\int_W \Psi(x_{2k})(\ddc\varphi_n+\Theta')
\wedge \tau^*(\Phi\wedge S)
\rightarrow \int_{\widetilde\Delta} \Phi\wedge S$$
for every test $(2k-p-1,2k-p-1)$-form $\Phi$ with compact support in 
$\widetilde\Delta\cap W$ and for every function $\Psi(x_{2k})$ 
supported in $\{|x_{2k}|<1\}$, such that $\Psi(0)=1$. 
Observe that since $\tau^*(\Phi\wedge S)$ is proportional to 
$\d x_1\wedge\d\overline x_1\wedge\ldots \wedge\d x_{2k-1} \wedge \d\overline x_{2k-1}$
only the component of $\ddc\varphi_n+\Theta'$ with respect to $\d x_{2k}\wedge 
\d \overline x_{2k}$ is relevant.
When $(x_1,\ldots, x_{2k-1})$
is fixed, we have
$$\int_{x_{2k}} \Psi(\d\d^{\rm c}_{x_{2k}}\varphi_n+\Theta')\rightarrow 1$$
since $\d\d^{\rm c}_{x_{2k}}\varphi_n+\Theta'$ converges to the Dirac mass
$\delta_0$ and $\Psi(0)=1$.
The last integral is uniformly bounded with respect
to $n$ and $x_1,\ldots,x_{2k-1}$ because by (6) one can prove that the masses of the measures 
$\d\d^{\rm c}_{x_{2k}}\varphi_n+\Theta'$
on a compact sets of $\{|x_{2k}|<1,\ x_1,\ldots , x_{2k-1}\mbox{ fixed}\}$ are uniformly bounded. This implies
the result.
\par
\hfill $\square$
\\
\begin{remark} \rm
Theorem 4.1 implies that on an arbitrary compact K\"ahler 
manifold $(X,\omega)$ 
positive pluriharmonic currents $T$ of bidegree $(1,1)$ have finite
energy. We then have $T=\Omega+\partial S +\overline\partial\overline S +
i\partial\overline\partial v$ with $\Omega$ smooth closed, 
$S$, $\partial S$ , $\overline\partial S$
in $\Ltwo$ and $v$ in $\Lone$. The energy of $T$ is equal to
$\int\overline\partial S\wedge \partial\overline S\wedge\omega^{k-2}$. 
The case of the projective space is treated in  \cite{FornaessSibony}.
To extend the result to an arbitrary compact K\"ahler manifold, one has to use the 
approximation Theorem 4.1, to go from a priori estimates on smooth positive pluriharmonic forms to
the estimates on positive pluriharmonic currents.
\end{remark}

Let $\DSH^p(X)$ denote the space of $(p,p)$-currents $T=T_1-T_2$ 
where $T_i$ are negative currents, such that 
$\ddc T_i=\Omega_i^+-\Omega_i^-$ with $\Omega_i^\pm$ 
positive closed. Observe that $\|\Omega_i^+\|=\|\Omega_i^-\|$. 
We define the {\it DSH-norm}  of $T$ as 
$$\|T\|_\DSH:=
\min\{\|T_1\|+\|T_2\|+\|\Omega_1^+\|+\|\Omega_2^+\|,\ T_i,\ \Omega_i^\pm
\mbox{ as above}\}.$$
We say that $T_n\rightarrow T$ in $\DSH^p(X)$ 
if $T_n\rightarrow T$ weakly and $(\|T_n\|_\DSH)$ is
bounded.

The spaces $\DSH^p(X)$ are analoguous to the space generated by q.p.s.h. functions.
They are useful in order to study the regularity of Green currents in dynamics \cite{DinhSibony3}. 
The proof of the following theorem, which gives the density of smooth
forms in $\DSH^p(X)$, follows the lines of previous approximation results and 
is left to the reader. In this case, for the control of the mass of $T_n^\pm$ we need to estimate the mass
of $\ddc\varphi_n\wedge\Pi^*(T)$. It is sufficient to estimate the mass of 
$\varphi_n \Pi^*(\ddc T)$ using the definition of $\varphi_n$. 
\begin{theorem} Let $T$ be a current 
in $\DSH^p(X)$. Then there exist smooth real $(p,p)$-forms 
$T_n$ such that $T_n\rightarrow T$. Moreover, $\|T_n\|_\DSH\leq c_X
\|T\|_\DSH$
where $c_X>0$ is a constant independent of $T$.  
\end{theorem}

\begin{remark}\rm
We have $\widetilde K_n^+ -\widetilde K_n^- -\gamma\wedge 
[\widetilde\Delta]= \gamma\wedge 
\ddc(\varphi_n-\varphi)$ and $\varphi_n=\varphi$ out 
of the set $(\varphi<-n+2)$. Hence
$\supp(\widetilde K_n^+ -\widetilde K_n^-)$ converge to $\widetilde \Delta$,
$\supp(K_n^+-K_n^-)$ converge to $\Delta$ and $\supp(T_n^+-T_n^-)$ 
converge to $\supp(T)$.
\end{remark}

We also have the following useful proposition.

\begin{proposition} Let $T$ be a continuous form and $T_n^\pm$
be the forms defined in (7)(8). Then $T_n:=T_n^+-T_n^-$ converge uniformly to $T$.
\end{proposition}
\begin{proof} We can approximate $T$ uniformly by smooth forms. We then assume 
that $T$ is smooth (see Lemma 3.1). 
The form $T_n-T$ is the push-forward of 
$(\widetilde K_n^+ -\widetilde K_n^- -\gamma\wedge 
[\widetilde\Delta])\wedge \Pi^*(T)$ by $\Pi':=\pi_1\circ\pi$.  
The last current is equal 
to $\widetilde T_n:=\ddc(\varphi_n-\varphi)\wedge \gamma'$ 
where $\gamma'$ is a smooth form.
Using a partition of unity, we reduce the problem to a local situation 
with the coordinates $x=(x',x'')=(x_1,\ldots,x_k,x_{k+1},\ldots, x_{2k})$, 
$\Pi'(x)=x'$, $\widetilde\Delta=(x_{2k}=0)$ and $\gamma'$ of compact support
as in the proof of Theorem 4.1. We have to check that 
$\Pi'_*(\widetilde T_n)(x')=\int_{x''} \widetilde T_n(x)$ converge uniformly to 0.

Observe that the last integral is taken in the neighbourhood 
$(\varphi<-n+2)$ of $(x_{2k}=0)$ and the form $\widetilde T_n - 
\d\d^{\rm c}_{x''}(\varphi_n-\varphi)\wedge \gamma'$ is of order $1/|x_{2k}|$
since in the difference we get at most one derivative with respect to $x_{2k}$. Hence,
it is sufficient to estimate 
$\int_{x''} \d\d^{\rm c}_{x''}(\varphi_n-\varphi)\wedge \gamma'= 
\int_{x''} (\varphi_n-\varphi)\wedge \d\d^{\rm c}_{x''}\gamma'$. 
It is clear that these forms
converge uniformly to 0.
\end{proof}

\section{Intersection of currents}

We want to consider a class of positive pluriharmonic currents
which are of interest in some problems of complex analysis and dynamics.
Some of their properties are 
given in \cite{Sibony2, AlessandriniBassanelli, 
Bassanelli, DuvalSibony, DinhLawrence, FornaessSibony}. Given a compact
K\"ahler manifold $(X,\omega)$ of dimension $k$, we want to define the 
intersection $S\wedge T$ of a positive closed $(1,1)$-current $S$
with a positive pluriharmonic current $T$ of bidegree 
$(p,p)$, $1\leq p\leq k-1$.
We have seen a special case of this situation in the last section.  

We write $S=\alpha+\ddc u$ with $\alpha$ smooth and $u$ a q.p.s.h. function. We
say that $u$ is {\it a potential} of $S$. 
\begin{theorem} Assume that
$u$ is continuous. Then $S\wedge T$ is well defined
and is a positive $\ddc$-closed current.
Moreover $S\wedge T$ depends continuously on $S$ and $T$ in the following sense. 
Let $T_n$ be positive pluriharmonic currents converging weakly to $T$. 
If $S_n=\alpha+\ddc u_n$ with $u_n$ continuous converging uniformly to 
$u$ then $S_n\wedge T_n$ converges weakly to $S\wedge T$. In particular, it holds
when the $u_n$ are continuous and decrease to $u$.
\end{theorem}
We first prove the following proposition for smooth potentials. We will see later that 
it can be 
extended to continuous q.p.s.h. functions $v^\pm$ and that 
$\d v^\pm\wedge R$ and $\dc v^\pm\wedge R$ are well defined in this case.
\begin{proposition} Let $v^\pm$ and $v=v^+-v^-$  be smooth real
functions on $X$ such 
that $\ddc v^\pm=\Theta^\pm-\alpha$ where $\alpha$ is a smooth closed $(1,1)$-form and 
$\Theta^\pm$ are positive closed $(1,1)$-currents. Let  
$R$ be a positive current in $\DSH^p(X)$ with $\ddc R=\Omega^+-\Omega^-$ where
$\Omega^\pm$ are positive closed currents. Then
\begin{eqnarray*}
\lefteqn{\int \d v\wedge\dc v\wedge R\wedge \omega^{k-p-1}
 \leq  }\\
&\leq &  
\|v\|_\Linfty \left(2\int \alpha\wedge R\wedge \omega^{k-p-1} + 
3(\|v^+\|_\Linfty  +\|v^-\|_\Linfty)\|\Omega^\pm\|\right).
\end{eqnarray*}
In particular, if $R$ is positive pluriharmonic, we have
$$\int \d v\wedge\dc v\wedge R\wedge \omega^{k-p-1}
 \leq 2\|v\|_\Linfty\int[\alpha]\wedge [R]\wedge [\omega]^{k-p-1}.$$ 
\end{proposition}
\begin{proof} Observe that $\|v\|_\Linfty\leq \|v^+\|_\Linfty +\|v^-\|_\Linfty$, 
$\|\Theta^+\|=\|\Theta^-\|$ and $\|\Omega^+\| = 
\|\Omega^-\|$. Hence
\begin{eqnarray*}
\lefteqn{\int \d v\wedge \dc v\wedge R\wedge\omega^{k-p-1}\leq}\\ 
& \leq & 
\frac{1}{2}\left|\int \ddc v^2\wedge R\wedge \omega^{k-p-1}\right| + 
\left|\int v\ddc v \wedge R\wedge \omega^{k-p-1} \right|\\
& = & \frac{1}{2}\left|\int v^2\wedge \ddc R\wedge \omega^{k-p-1}\right|
+ \left|\int v(\Theta^+-\Theta^-)\wedge R\wedge \omega^{k-p-1}\right| \\
& \leq & \frac{1}{2}\|v\|_\Linfty^2 \int (\Omega^++\Omega^-)\wedge \omega^{k-p-1}
+\|v\|_\Linfty \int (\Theta^++\Theta^-)\wedge R\wedge \omega^{k-p-1} \\
& = & \|v\|_\Linfty^2 \|\Omega^\pm\| + 
2\|v\|_\Linfty\int \alpha\wedge R\wedge \omega^{k-p-1}+\\
& & 
+ \|v\|_\Linfty\int \ddc (v^+ +v^-) \wedge R\wedge \omega^{k-p-1}\\
& \leq & \|v\|_\Linfty ( \|v^+\|_\Linfty  +\|v^-\|_\Linfty)\|\Omega^\pm\| + 
2\|v\|_\Linfty\int \alpha\wedge R\wedge \omega^{k-p-1}+\\
& & 
+ \|v\|_\Linfty\int (v^+ +v^-) \wedge (\Omega^+-\Omega^-)\wedge \omega^{k-p-1}\\
&\leq & \|v\|_\Linfty \left(2\int \alpha\wedge R\wedge \omega^{k-p-1} + 
3(\|v^+\|_\Linfty  +\|v^-\|_\Linfty)\|\Omega^\pm\|\right).
\end{eqnarray*}
\end{proof}
{\it Proof of Theorem 5.1.}
Observe that, when $u_n$ decreases to $u$, the Hartogs lemma implies that
$u_n$ converges uniformly to $u$. 
By Demailly's regularization theorem \cite{Demailly2}, we can assume that
$u_n$ are smooth and uniformly convergent to $u$. 
So, $S_n\wedge T_n$ is well defined. We will prove that 
$S_n\wedge T_n$ converges. This also implies that the limit depends 
only on $S$ and $T$.

We first consider the case where $T_n=T$. We then have
\begin{eqnarray}
S_n\wedge T = \alpha\wedge T +\d(\dc u_n\wedge T)-\dc(\d u_n\wedge T) 
- \ddc (u_n T).
\end{eqnarray}
Proposition 5.2 applied to $u_n-u_m$ and the Cauchy criterion imply that 
$\d u_n$ and $\dc u_n$ converge in 
$\Ltwo(T\wedge \omega^{k-p-1})$. Hence $S_n\wedge T$ converges and
$\d u\wedge T$, $\dc u\wedge T$ are well defined.

To complete the proof we write 
$$S_n\wedge T_n-S \wedge T = \ddc (u_n-u)\wedge T_n + 
\ddc u\wedge (T_n-T) +\alpha\wedge(T_n-T).$$
The last term tends to zero. 
Proposition 5.2 implies that $\int \d(u_n-u) \wedge \dc(u_n-u) \wedge T_n
\wedge\omega^{k-p-1}$ has zero limit.
An identity as in (12) shows that the first term tends to $0$. 
For the second term, we observe
first that if $\gamma$ is a test $1$-form and $v$ is a smooth q.p.s.h. function
with $\|u-v\|_\Linfty \leq \epsilon$, then
\begin{eqnarray*}
\int\d u\wedge\gamma\wedge(T-T_n)\wedge \omega^{k-p-1}
& = & \int\d v\wedge\gamma\wedge(T-T_n)\wedge \omega^{k-p-1}+\\
& & 
+ \int\d (u-v)\wedge\gamma\wedge(T-T_n)\wedge \omega^{k-p-1}.
\end{eqnarray*}
The first integral tends to zero. Schwarz's inequality and Proposition 5.2 
imply
\begin{eqnarray*}
\lefteqn{\left|\int\d (u-v)\wedge\gamma\wedge T_n\wedge \omega^{k-p-1}\right|^2
 \leq }\\ 
&\leq & \mbox{const} \int \d (u-v) \wedge \dc (u-v) 
\wedge T_n\wedge \omega^{k-p-1}\\
& \leq &  \mbox{const} \|u-v\|_\Linfty\|T_n\| 
\end{eqnarray*}
and similarly for $T$.
\par
\hfill $\square$
\\

In the same way, using the full strength of Proposition 5.2, one 
can prove the following theorem.
\begin{theorem} Let $T$ be a current in $\DSH^p(X)$ and $S$ as in Theorem 5.1.
Then $S\wedge T$ is well defined
and belongs to $\DSH^{p+1}(X)$. Moreover $S\wedge T$ depends continuously on $S$ 
and $T$. The topology on the $T$ variable is the topology of $\DSH^p(X)$.
\end{theorem}
It is enough to assume $T$ positive and to modify (12) into
\begin{eqnarray*}
S_n\wedge T = \alpha\wedge T +\d(\dc u_n\wedge T)-\dc(\d u_n\wedge T) 
- \ddc (u_n T) +u_n\ddc T.
\end{eqnarray*}
\begin{remarks} \rm
If $S_i$ are positive closed $(1,1)$-currents with continuous potentials, then 
$S_1\wedge\ldots\wedge S_m\wedge T$ is symmetric in $S_i$ since 
this is true when $S_i$ and $T$ are smooth.
Let $u$ be a p.s.h. function on an open ball $\Omega\subset X$. By the maximum 
regularization procedure as in \cite{Demailly1}, 
if $u$ is continuous we can extend $u$ to a continuous 
q.p.s.h. function on $X$. 
Hence, $\ddc u\wedge T$ is well defined on $\Omega$. 

When $T$ is only a (positive pluriharmonic) $(p,p)$-current on $\Omega$, 
we don't know how to define $\ddc u\wedge T$ without additional hypothesis on $u$.
Assume that $T$, $\d T$
and $\ddc T$ are of order 0 and $u$ is
locally integrable with respect to the coefficient measures of 
$T$, $\d T$ and $\ddc T$. Then we can define
$$\ddc u\wedge T:=\ddc(uT)+u\ddc T -\d(u\dc T)+\dc(u\d T).$$
If $u_n$ are p.s.h. and decrease to $u$ or if 
$u_n$ converge uniformly to $u$, we have
$\ddc u_n\wedge T\rightarrow \ddc u\wedge T$. When $T$ is positive, we also
have an inequality of Chern-Levine-Nirenberg type (see \cite[p.126]{Demailly1}).
More precisely, if $K$, $L$ are compact sets in $\Omega$ with $L\Subset K$, 
then
$$\|\ddc u\wedge T\|_L\leq c_{K,L}(\|uT\|_K +\|u\d T\|_K +
\|u\ddc T\|_K).$$
Note that positive harmonic currents associated to a lamination by Riemann surfaces 
satisfy the above hypothesis (see \cite{BerndtssonSibony}).
\end{remarks}

If $T$ is of bidegree $(1,1)$ we can extend Theorem 5.1 to currents $S$ with bounded
potential.
\begin{proposition} Let $T$ be a positive pluriharmonic current of bidegree $(1,1)$
in $(X,\omega)$. If $u$ is a bounded q.p.s.h. function then $\ddc u\wedge T$ is well
defined. If $(u_n)$ is a bounded sequence of q.p.s.h. functions converging pointwise
to $u$ with $\ddc u_n\geq -c\omega$, then $\ddc u_n\wedge T\rightarrow \ddc u\wedge T$.  
\end{proposition}
\begin{proof} We can assume that $u_n$ are smooth and positive. It is easy to check 
that $u_n^2$ are q.p.s.h. and converge to $u^2$. It follows that 
$\partial u_n$ (resp. $\overline \partial u_n$) converge to 
$\partial u$ (resp. $\overline \partial u$) weakly in $\Ltwo(X)$.

As in Theorem 5.1, we only need to show that $\partial u_n\wedge T$ 
(resp. $\overline \partial u_n\wedge T$) converges weakly. Recall that we can write
$T=\Omega+\partial S+ \overline\partial\overline S + i\partial\overline \partial v$
with $\Omega$ smooth closed, 
$\partial S$, $\overline\partial  S$ 
in $\Ltwo$, $v$ in $\Lone$ \cite{FornaessSibony}.
If $\gamma$ is a test $1$-form, we have
\begin{eqnarray*}
\int\partial u_n\wedge T\wedge\gamma \wedge \omega^{k-1}
&  = &   - \int u_n \partial T\wedge \gamma \wedge \omega^{k-1}- 
\int u_n T\wedge \partial\gamma\wedge \omega^{k-1}\\
& = &  - \int u_n \partial \overline \partial \overline S 
\wedge \gamma \wedge \omega^{k-1}- 
\int u_n T\wedge \partial\gamma\wedge \omega^{k-1} 
\end{eqnarray*}
The second term tends to $\int u T\wedge \partial\gamma\wedge \omega^{k-1}$.
The first term is equal to
$$ - \int \partial u_n \wedge \overline \partial \overline S 
\wedge \gamma \wedge \omega^{k-1}+ 
\int u_n \overline \partial \overline S \wedge \partial\gamma\wedge \omega^{k-1}$$
which converges to
$$ - \int \partial u \wedge \overline \partial \overline S 
\wedge \gamma \wedge \omega^{k-1}+ 
\int u \overline \partial \overline S \wedge \partial\gamma\wedge \omega^{k-1}$$
since $u_n\rightarrow u$ and $\partial u_n\rightarrow\partial u$ weakly in $\Ltwo$.

The convergence of $\overline \partial u_n\wedge T$ is proved in the same way.
\end{proof}

\small

Tien-Cuong Dinh and Nessim Sibony,\\
Math\'ematique - B\^at. 425, UMR 8628, 
Universit\'e Paris-Sud, 91405 Orsay, France. \\
E-mails: TienCuong.Dinh@math.u-psud.fr and
Nessim.Sibony@math.u-psud.fr.

\begin{thebibliography}{11}
%
\bibitem{AlessandriniBassanelli}
\textit{L. Alessandrini and G. Bassanelli},
Plurisubharmonic currents and their extension across analytic subsets, 
\textit{Forum Math.}, \textbf{5} (1993), no. 6, 577--602.
%
\bibitem{Bassanelli}
\textit{G. Bassanelli}, A cut-off theorem for plurisubharmonic currents,  
\textit{Forum Math.},  \textbf{6}  (1994),  no. 5, 567--595.
%
\bibitem{BerndtssonSibony}
\textit{B. Berndtsson and N. Sibony},
The $\overline\partial$ equation
on a positive current, \textit{Invent. Math.},
\textbf{147} (2002), 371-428.
%
\bibitem{Blanchard}
\textit{A. Blanchard},
Sur les vari\'et\'es analytiques complexes, 
\textit{Ann. Sci. Ecole Norm. Sup.} (3),
\textbf{73} (1956), 157--202.
%
\bibitem{Bowen}
\textit{R. Bowen}, Topological entropy for non compact sets,
\textit{Trans. A.M.S.}, \textbf{184} (1973), 125-136.
%
\bibitem{Demailly1}
\textit{J.P. Demailly}, Monge-Amp\`ere Operators, Lelong numbers and
Intersection theory in Complex Analysis and Geometry, \textit{Plemum
  Press} (1993), 115-193, \textit{(V. Ancona and A. Silva editors)}.
%
\bibitem{Demailly2}
\textit{J.P. Demailly}, Pseudoconvex-concave duality and regularization of
currents. Several complex variables (Berkeley, CA, 1995-1996), 233-271, 
\textit{Math. Sci. Res. Inst. Publ.}, \textbf{37}, Cambridge Univ. Press, 
Cambridge, 1999. 
%
\bibitem{Dinh}
\textit{T.C. Dinh}, Distribution des pr\'eimages et des points
p\'eriodiques d'une correspondance polynomiale,
{\it Bull. Soc. Math. France}, to appear. 
%
\bibitem{DinhLawrence}
\textit{T.C. Dinh and M. Lawrence}, Polynomial hull and postive currents,
{\it Ann. Fac. Sci. Toulouse}, Vol. {\bf XII}, no 3 (2003), 317-334.
%
\bibitem{DinhSibony1}
\textit{T.C. Dinh and N. Sibony}, 
Une borne sup\'erieure pour l'entropie topologique d'une 
application rationnelle,
\textit{Ann. of Math.}, to appear. 
%
\bibitem{DinhSibony2}
\textit{T.C. Dinh and N. Sibony}, Distribution des valeurs de
transformations m\'eromorphes et applications,
\textit{preprint}, 2003. 
arxiv.org/abs/math.DS/0306095. 
%
\bibitem{DinhSibony3}
\textit{T.C. Dinh and N. Sibony}, Green currents for automorphisms of 
compact K\"ahler manifolds, \textit{JAMS},
to appear.
%
\bibitem{DuvalSibony}
{\it J. Duval and N. Sibony}, 
Polynomial convexity, rational convexity, and currents,
\textit{Duke Math. J.}, \textbf{79}, No.2 (1995), 487-513. 
%
\bibitem{FornaessSibony}
\textit{J.E. Forn\ae ss and N. Sibony}, Harmonic currents with finite energy,
\textit{preprint}, 2004. arxiv.org/abs/math.CV/0402432.
%
\bibitem{Friedland}
\textit{S. Friedland}, Entropy of polynomial and rational maps,
\textit{Ann. of Math.}, \textbf{133} (1991), 359-368.
%
\bibitem{Gromov1}
\textit{M. Gromov}, On the entropy of holomorphic maps,
\textit{S\'eminaire Bourbaki}, vol. 1985/86, \textit{Ast\'erique}, 
vol. {\bf 145-146}, Soc. Math. France, 1987, 225-240.
%
\bibitem{Gromov2}
\textit{M. Gromov}, On the entropy of holomorphic maps,
\textit{Enseignement Math.}, \textbf{49} (2003), 217-235. {\it Manuscript} (1977). 
%
\bibitem{Lelong}
\textit{P. Lelong}, Fonctions plurisousharmoniques et formes
diff\'erentielles positives, Dunod Paris, 1968.
%
\bibitem{Meo}
\textit{M. M\'eo}, Image inverse d'un courant positif ferm\'e par une
application surjective, \textit{C.R.A.S.}, \textbf{322} (1996),
1141-1144. 
%
\bibitem{Sibony}
\textit{N. Sibony}, Dynamique des applications rationnelles de
$\mathbb{P}^k$, \textit{Panoramas et Synth\`eses}, {\bf 8} (1999), 97-185.
%
\bibitem{Sibony2}
\textit{N. Sibony}, Quelques probl\`emes de 
prolongement de courants en analyse complexe, 
\textit{Duke Math. J.}, \textbf{52} (1985),  no.1, 157--197.
%
\bibitem{Skoda}
\textit{H. Skoda}, Prolongement des courants positifs, ferm\'es de
masse finie, \textit{Invent. Math.}, \textbf{66} (1982), 361-376.
%
\bibitem{Yomdin}
\textit{Y. Yomdin}, Volume growth and entropy, \textit{Israel
  J. Math.}, \textbf{57} (1987), 285-300.
%
\end{thebibliography}
\end{document}